\documentclass{amsart}%
\usepackage{amsfonts}
\usepackage{graphicx}
\usepackage{amscd}%
\usepackage{amsmath}%
\usepackage{amssymb}
\newtheorem{theorem}{Theorem}
\theoremstyle{plain}

\newtheorem{corollary}{Corollary}

\newtheorem{lemma}{Lemma}

\numberwithin{equation}{section}
\numberwithin{theorem}{section}
\numberwithin{algorithm}{section}
\numberwithin{axiom}{section}
\numberwithin{case}{section}
\numberwithin{claim}{section}
\numberwithin{conclusion}{section}
\numberwithin{condition}{section}
\numberwithin{conjecture}{section}
\numberwithin{corollary}{section}
\numberwithin{criterion}{section}
\numberwithin{definition}{section}
\numberwithin{example}{section}
\numberwithin{exercise}{section}
\numberwithin{lemma}{section}
\numberwithin{notation}{section}
\numberwithin{problem}{section}
\numberwithin{proposition}{section}
\numberwithin{remark}{section}
\numberwithin{solution}{section}

\begin{document}
\title[Weak limits of smooth maps]{On the weak limits of smooth maps for the Dirichlet energy between manifolds}
\author{Fengbo Hang}
\address{Department of Mathematics, Princeton University, Fine Hall, Washington Road,
Princeton, NJ 08544, and, School of Mathematics, Institute for Advanced Study,
1 Einstein Drive, Princeton, NJ 08540}
\email{fhang@math.princeton.edu}

\begin{abstract}
We identify all the weak sequential limits of smooth maps in $W^{1,2}\left(
M,N\right)  $. In particular, this implies a necessary and sufficient
topological condition for smooth maps to be weakly sequentially dense in
$W^{1,2}\left(  M,N\right)  $.

\end{abstract}
\maketitle

\section{Introduction\label{sec1}}

Assume $M$ and $N$ are smooth compact Riemannian manifolds without boundary
and they are embedded into $\mathbb{R}^{l}$ and $\mathbb{R}^{\overline{l}}$
respectively. The following spaces are of interest in the calculus of
variations:%
\begin{align*}
W^{1,2}\left(  M,N\right)   & =\left\{  u\in W^{1,2}\left(  M,\mathbb{R}%
^{\overline{l}}\right)  :u\left(  x\right)  \in N\text{ a.e. }x\in M\right\}
,\\
H_{W}^{1,2}\left(  M,N\right)   & =\left\{  u\in W^{1,2}\left(  M,N\right)
:\text{ there exists a sequence }u_{i}\in C^{\infty}\left(  M,N\right)
\right. \\
& \left.  \text{such that }u_{i}\rightharpoonup u\text{ in }W^{1,2}\left(
M,N\right)  \right\}  .
\end{align*}
For a brief history and detailed references on the study of analytical and
topological issues related to these spaces, one may refer to \cite{HL1,HL2,PR}%
. In particular, it follows from theorem 7.1 of \cite{HL2} that a necessary
condition for $H_{W}^{1,2}\left(  M,N\right)  =W^{1,2}\left(  M,N\right)  $ is
that $M$ satisfies the $1$-extension property with respect to $N$ (see section
2.2 of \cite{HL2} for a definition). It was conjectured in section 7 of
\cite{HL2} that the $1$-extension property is also sufficient for $H_{W}%
^{1,2}\left(  M,N\right)  =W^{1,2}\left(  M,N\right)  $. In \cite{H,PR}, it
was shown that $H_{W}^{1,2}\left(  M,N\right)  =W^{1,2}\left(  M,N\right)  $
when $\pi_{1}\left(  M\right)  =0$ or $\pi_{1}\left(  N\right)  =0 $. Note
that if $\pi_{1}\left(  M\right)  =0$ or $\pi_{1}\left(  N\right)  =0$, then
$M$ satisfies the $1$-extension property with respect to $N$. In section 8 of
\cite{HL3}, it was proved that the above conjecture is true under the
additional assumption that $N$ satisfies the $2$-vanishing condition. The main
aim of the present article is to confirm the conjecture in its full
generality. More precisely, we have

\begin{theorem}
\label{thm1.1}Let $M^{n}$ and $N$ be smooth compact Riemannian manifolds
without boundary ($n\geq3$). Take a Lipschitz triangulation $h:K\rightarrow
M$, then%
\begin{align*}
& H_{W}^{1,2}\left(  M,N\right) \\
& =\left\{  u\in W^{1,2}\left(  M,N\right)  :u_{\#,2}\left(  h\right)  \text{
has a continuous extension to }M\text{ w.r.t. }N\right\} \\
& =\left\{  u\in W^{1,2}\left(  M,N\right)  :u\text{ may be connected to some
smooth maps}\right\}  .
\end{align*}
In addition, if $\alpha\in\left[  M,N\right]  $ satisfies $\left.  \alpha\circ
h\right\vert _{\left\vert K^{1}\right\vert }=u_{\#,2}\left(  h\right)  $, then
we may find a sequence of smooth maps $u_{i}\in C^{\infty}\left(  M,N\right)
$ such that $u_{i}\rightharpoonup u$ in $W^{1,2}\left(  M,N\right)  $,
$\left[  u_{i}\right]  =\alpha$ and $du_{i}\rightarrow du$ a.e..
\end{theorem}

Here $u_{\#,2}\left(  h\right)  $ is the $1$-homotopy class defined by White
\cite{W} (see also section 4 of \cite{HL2}) and $\left[  M,N\right]  $ means
all homotopy classes of maps from $M$ to $N$. It follows from Theorem
\ref{thm1.1} that

\begin{corollary}
\label{cor1.1}Let $M^{n}$ and $N$ be smooth compact Riemannian manifolds
without boundary and $n\geq3$. Then smooth maps are weakly sequentially dense
in $W^{1,2}\left(  M,N\right)  $ if and only if $M$ satisfies the
$1$-extension property with respect to $N$.
\end{corollary}

For $p\in\left[  3,n-1\right]  $ being an natural number, it remains a
challenging open problem to find out whether the weak sequential density of
smooth maps in $W^{1,p}\left(  M,N\right)  $ is equivalent to the condition
that $M$ satisfies the $p-1$ extension property with respect to $N$. This was
verified to be true under further topological assumptions on $N$ (see section
8 of \cite{HL3}). However, even for $W^{1,3}\left(  S^{4},S^{2}\right)  $, it
is still not known whether smooth maps are weakly sequentially dense. Some
very interesting recent work on this space can be found in \cite{HR}.

The paper is written as follows. In Section \ref{sec2}, we will present some
technical lemmas. In Section \ref{sec3}, we will prove the above theorem and corollary.

\textbf{Acknowledgments. }The research of the author is supported by National
Science Foundation Grant DMS-0209504.

\section{Some preparations\label{sec2}}

The following local result, which was proved by Pakzad and Riviere in
\cite{PR}, plays an important role in our discussion.

\begin{theorem}
[\cite{PR}]\label{thmPR}Let $N$ be a smooth compact Riemannian manifold.
Assume $n\geq3 $, $B_{1}=B_{1}^{n}$, $f\in W^{1,2}\left(  \partial
B_{1},N\right)  \cap C\left(  \partial B_{1},N\right)  $, $f\sim
\operatorname{const}$, $u\in W^{1,2}\left(  B_{1},N\right)  $, $\left.
u\right\vert _{\partial B_{1}}=f$, then there exists a sequence $u_{i}\in
W^{1,2}\left(  B_{1},N\right)  \cap C\left(  \overline{B}_{1},N\right)  $ such
that $\left.  u_{i}\right\vert _{\partial B_{1}}=f$, $u_{i}\rightharpoonup u$
in $W^{1,2}\left(  B_{1},N\right)  $ and $du_{i}\rightarrow du$ a.e.. In
addition, if $v\in W^{1,2}\left(  B_{2}\backslash B_{1},N\right)  \cap
C\left(  \overline{B}_{2}\backslash B_{1},N\right)  $ satisfies $\left.
v\right\vert _{\partial B_{1}}=f$ and $\left.  v\right\vert _{\partial B_{2}%
}\equiv\operatorname{const}$, then we may estimate%
\[
\int_{B_{1}}\left\vert du_{i}\right\vert ^{2}d\mathcal{H}^{n}\leq c\left(
n,N\right)  \left(  \int_{B_{1}}\left\vert du\right\vert ^{2}d\mathcal{H}%
^{n}+\int_{B_{2}\backslash B_{1}}\left\vert dv\right\vert ^{2}d\mathcal{H}%
^{n}\right)  .
\]

\end{theorem}

For convenience, we will use those notations and concepts in section 2, 3 and
4 of \cite{HL2}. The following lemma is a rough version of Luckhaus's lemma
\cite{L}. For reader's convenience, we sketch a proof of this simpler version
using results from section 3 of \cite{HL2}.

\begin{lemma}
\label{lem2.1}Assume $M^{n}$ and $N$ are smooth compact Riemannian manifolds
without boundary. Let $e>0$, $0<\delta<1$, $A>0$, then there exists an
$\varepsilon=\varepsilon\left(  e,\delta,A,M,N\right)  >0$ such that for any
$u,v\in W^{1,2}\left(  M,N\right)  $ with $\left\vert du\right\vert
_{L^{2}\left(  M\right)  },\left\vert dv\right\vert _{L^{2}\left(  M\right)
}\leq A$ and $\left\vert u-v\right\vert _{L^{2}\left(  M\right)  }%
\leq\varepsilon$, we may find a $w\in W^{1,2}\left(  M\times\left(
0,\delta\right)  ,N\right)  $ such that, in the trace sense $w\left(
x,0\right)  =u\left(  x\right)  $, $w\left(  x,\delta\right)  =v\left(
x\right)  $ a.e. $x\in M$ and%
\[
\left\vert dw\right\vert _{L^{2}\left(  M\times\left(  0,\delta\right)
\right)  }\leq c\left(  M\right)  \sqrt{\delta}\left(  \left\vert
du\right\vert _{L^{2}\left(  M\right)  }+\left\vert dv\right\vert
_{L^{2}\left(  M\right)  }+e\right)  .
\]

\end{lemma}

\begin{proof}
Let $\varepsilon_{M}>0$ be a small positive number such that%
\[
V_{2\varepsilon_{M}}\left(  M\right)  =\left\{  x\in\mathbb{R}^{l}:d\left(
x,M\right)  <2\varepsilon_{M}\right\}
\]
is a tubular neighborhood of $M$. Let $\pi_{M}:V_{2\varepsilon_{M}}\left(
M\right)  \rightarrow M$ be the nearest point projection. Similarly we have
$\varepsilon_{N}$, $V_{2\varepsilon_{N}}\left(  N\right)  $ and $\pi_{N}$ for
$N$. Choose a Lipschitz cubeulation $h:K\rightarrow M$. We may assume each
cell in $K$ is a cube of unit size. For $\xi\in B_{\varepsilon_{M}}^{l}$,
$x\in\left\vert K\right\vert $, let $h_{\xi}\left(  x\right)  =\pi_{M}\left(
h\left(  x\right)  +\xi\right)  $. Assume $\varepsilon_{M}$ is small enough
such that all $h_{\xi}$'s are bi-Lipschitz maps. Set $m=\left[  \frac
{1}{\delta}\right]  +1$, using $\left[  0,1\right]  =\cup_{i=1}^{m}\left[
\frac{i-1}{m},\frac{i}{m}\right]  $, we may divide each $k$-cube in $K$ into
$m^{k}$ small cubes. In particular, we get a subdivision of $K$, called
$K_{m}$. It follows from section 3 of \cite{HL2} that for a.e. $\xi\in
B_{\varepsilon_{M}}^{l}$, $u\circ h_{\xi},v\circ h_{\xi}\in\mathcal{W}%
^{1,2}\left(  K_{m},N\right)  $. Applying the estimates in section 3 of
\cite{HL2} to each unit size $k$-cube in $\left\vert K_{m}^{k}\right\vert $,
we get%
\begin{align*}
\int_{B_{\varepsilon_{M}}^{l}}d\mathcal{H}^{l}\left(  \xi\right)
\int_{\left\vert K_{m}^{k}\right\vert }\left\vert d\left(  \left.  u\circ
h_{\xi}\right\vert _{\left\vert K_{m}^{k}\right\vert }\right)  \right\vert
^{2}d\mathcal{H}^{k}  & \leq c\left(  M\right)  \delta^{k-n}\left\vert
du\right\vert _{L^{2}\left(  M\right)  }^{2},\\
\int_{B_{\varepsilon_{M}}^{l}}d\mathcal{H}^{l}\left(  \xi\right)
\int_{\left\vert K_{m}^{k}\right\vert }\left\vert d\left(  \left.  v\circ
h_{\xi}\right\vert _{\left\vert K_{m}^{k}\right\vert }\right)  \right\vert
^{2}d\mathcal{H}^{k}  & \leq c\left(  M\right)  \delta^{k-n}\left\vert
dv\right\vert _{L^{2}\left(  M\right)  }^{2},
\end{align*}
and%
\begin{align*}
& \left(  \int_{B_{\varepsilon_{M}}^{l}}\left\vert u\circ h_{\xi}-v\circ
h_{\xi}\right\vert _{L^{\infty}\left(  \left\vert K_{m}^{1}\right\vert
\right)  }^{2}d\mathcal{H}^{l}\left(  \xi\right)  \right)  ^{\frac{1}{2}}\\
& \leq c\left(  \delta,M\right)  \left(  \left\vert d\left(  u-v\right)
\right\vert _{L^{2}\left(  M\right)  }^{\frac{3}{4}}\left\vert u-v\right\vert
_{L^{2}\left(  M\right)  }^{\frac{1}{4}}+\left\vert u-v\right\vert
_{L^{2}\left(  M\right)  }\right) \\
& \leq c\left(  \delta,A,M\right)  \varepsilon^{\frac{1}{4}}.
\end{align*}
By the mean value inequality, we may find a $\xi\in B_{\varepsilon_{M}}^{l} $
such that $u\circ h_{\xi},v\circ h_{\xi}\in\mathcal{W}^{1,2}\left(
K_{m},N\right)  $,%
\[
\left\vert u\circ h_{\xi}-v\circ h_{\xi}\right\vert _{L^{\infty}\left(
\left\vert K_{m}^{1}\right\vert \right)  }\leq c\left(  \delta,A,M\right)
\varepsilon^{\frac{1}{4}}<\varepsilon_{N}\quad\text{when }\varepsilon\text{ is
small enough,}%
\]
and%
\begin{align*}
& \int_{\left\vert K_{m}^{k}\right\vert }\left[  \left\vert d\left(  \left.
u\circ h_{\xi}\right\vert _{\left\vert K_{m}^{k}\right\vert }\right)
\right\vert ^{2}+\left\vert d\left(  \left.  v\circ h_{\xi}\right\vert
_{\left\vert K_{m}^{k}\right\vert }\right)  \right\vert ^{2}\right]
d\mathcal{H}^{k}\\
& \leq c\left(  M\right)  \delta^{k-n}\left(  \left\vert du\right\vert
_{L^{2}\left(  M\right)  }^{2}+\left\vert dv\right\vert _{L^{2}\left(
M\right)  }^{2}\right)
\end{align*}
for $1\leq k\leq n$. Fix a $\eta\in C^{\infty}\left(  \mathbb{R}%
,\mathbb{R}\right)  $ such that $0\leq\eta\leq1$, $\left.  \eta\right\vert
_{\left(  -\infty,\frac{1}{3}\right)  }=1$ and $\left.  \eta\right\vert
_{\left(  \frac{2}{3},\infty\right)  }=0$. Letting $f=u\circ h_{\xi}$,
$g=v\circ h_{\xi}$, we will define $\phi:\left\vert K\right\vert \times\left[
0,\delta\right]  \rightarrow N$ inductively. First set $\phi\left(
x,0\right)  =f\left(  x\right)  $ and $\phi\left(  x,\delta\right)  =g\left(
x\right)  $ for $x\in\left\vert K\right\vert $. For $\Delta\in K_{m}%
^{1}\backslash K_{m}^{0}$, on $\Delta\times\left[  0,\delta\right]  $, we let
\[
\phi\left(  x,t\right)  =\pi_{N}\left(  \eta\left(  \frac{t}{\delta}\right)
f\left(  x\right)  +\left(  1-\eta\left(  \frac{t}{\delta}\right)  \right)
g\left(  x\right)  \right)  \quad x\in\Delta,0\leq t\leq\delta.
\]
For $\Delta\in K_{m}^{2}\backslash K_{m}^{1}$, let $y_{\Delta}$ be the center
of $\Delta$, and define $\phi$ on $\Delta\times\left[  0,\delta\right]  $ as
the homogeneous degree zero extension of $\left.  \phi\right\vert
_{\partial\left(  \Delta\times\left[  0,\delta\right]  \right)  } $ with
respect to $\left(  y_{\Delta},\frac{\delta}{2}\right)  $. Next we handle each
$3$-cube, $4$-cube, $\cdots$, $n$-cube in a similar way. Calculations show
that%
\begin{align*}
& \int_{\left\vert K\right\vert \times\left[  0,\delta\right]  }\left\vert
d\phi\right\vert ^{2}d\mathcal{H}^{n+1}\\
& \leq c\left(  n\right)  \sum_{k=1}^{n}\delta^{n+1-k}\int_{\left\vert
K_{m}^{k}\right\vert }\left[  \left\vert d\left(  \left.  u\circ h_{\xi
}\right\vert _{\left\vert K_{m}^{k}\right\vert }\right)  \right\vert
^{2}+\left\vert d\left(  \left.  v\circ h_{\xi}\right\vert _{\left\vert
K_{m}^{k}\right\vert }\right)  \right\vert ^{2}\right]  d\mathcal{H}%
^{k}+c\left(  \delta,A,M\right)  \varepsilon^{\frac{1}{2}}\\
& \leq c\left(  M\right)  \delta\left(  \left\vert du\right\vert
_{L^{2}\left(  M\right)  }^{2}+\left\vert dv\right\vert _{L^{2}\left(
M\right)  }^{2}+e^{2}\right)
\end{align*}
when $\varepsilon$ is small enough. Finally $w:M\times\left[  0,\delta\right]
\rightarrow N$, defined by $w\left(  x,t\right)  =\phi\left(  h_{\xi}%
^{-1}\left(  x\right)  ,t\right)  $, is the needed map.
\end{proof}

\begin{lemma}
\label{lem2.2}Assume $N$ is a smooth compact Riemannian manifold, $n\geq2$,
$B_{1}=B_{1}^{n}$, $u,v\in W^{1,2}\left(  B_{1},N\right)  $ such that $\left.
u\right\vert _{\partial B_{1}}=\left.  v\right\vert _{\partial B_{1}}$. Define
$w:B_{1}\times\left(  0,1\right)  \rightarrow N$ by%
\[
w\left(  x,t\right)  =\left\{
\begin{array}
[c]{cc}%
u\left(  x\right)  , & x\in B_{1}\backslash B_{t};\\
u\left(  \frac{t^{2}}{\left\vert x\right\vert }\frac{x}{\left\vert
x\right\vert }\right)  , & x\in B_{t}\backslash B_{t^{2}};\\
v\left(  \frac{x}{t^{2}}\right)  , & x\in B_{t^{2}};
\end{array}
\right.
\]
then $w\in W^{1,2}\left(  B_{1}\times\left(  0,1\right)  ,N\right)  $ and%
\[
\left\vert dw\right\vert _{L^{2}\left(  B_{1}\times\left(  0,1\right)
\right)  }\leq c\left(  n\right)  \left(  \left\vert du\right\vert
_{L^{2}\left(  B_{1}\right)  }+\left\vert dv\right\vert _{L^{2}\left(
B_{1}\right)  }\right)  .
\]

\end{lemma}

\begin{proof}
Note that%
\[
\left\vert dw\left(  x,t\right)  \right\vert \leq\left\{
\begin{array}
[c]{cc}%
\left\vert du\left(  x\right)  \right\vert , & t<\left\vert x\right\vert ;\\
c\left(  n\right)  \left\vert du\left(  \frac{t^{2}}{\left\vert x\right\vert
}\frac{x}{\left\vert x\right\vert }\right)  \right\vert \frac{t^{2}%
}{\left\vert x\right\vert ^{2}}, & t^{2}<\left\vert x\right\vert <t;\\
c\left(  n\right)  \left\vert dv\left(  \frac{x}{t^{2}}\right)  \right\vert
\frac{1}{t^{2}}, & \left\vert x\right\vert <t^{2}.
\end{array}
\right.
\]
Hence%
\begin{align*}
& \int_{\substack{0<t<1 \\t^{2}<\left\vert x\right\vert <t}}\left\vert
dw\left(  x,t\right)  \right\vert ^{2}d\mathcal{H}^{n+1}\left(  x,t\right) \\
& \leq c\left(  n\right)  \int_{0}^{1}dt\int_{t^{2}}^{t}dr\int_{\partial
B_{r}}\left\vert du\left(  \frac{t^{2}}{r^{2}}x\right)  \right\vert ^{2}%
\frac{t^{4}}{r^{4}}d\mathcal{H}^{n-1}\left(  x\right) \\
& =c\left(  n\right)  \int_{0}^{1}dt\int_{t}^{1}ds\int_{\partial B_{s}}%
\frac{t^{2\left(  n-2\right)  }}{s^{2\left(  n-2\right)  }}\left\vert
du\left(  y\right)  \right\vert ^{2}d\mathcal{H}^{n-1}\left(  y\right) \\
& \leq c\left(  n\right)  \left\vert du\right\vert _{L^{2}\left(
B_{1}\right)  }^{2},
\end{align*}
and%
\begin{align*}
& \int_{\substack{0<t<1 \\\left\vert x\right\vert <t^{2}}}\left\vert dw\left(
x,t\right)  \right\vert ^{2}d\mathcal{H}^{n+1}\left(  x,t\right) \\
& \leq c\left(  n\right)  \int_{0}^{1}dt\int_{B_{t^{2}}}\left\vert dv\left(
\frac{x}{t^{2}}\right)  \right\vert ^{2}\frac{1}{t^{4}}d\mathcal{H}^{n}\left(
x\right) \\
& \leq c\left(  n\right)  \left\vert dv\right\vert _{L^{2}\left(
B_{1}\right)  }^{2}.
\end{align*}
The lemma follows.
\end{proof}

\section{Identifying weak limits of smooth maps\label{sec3}}

In this section, we shall prove Theorem \ref{thm1.1} and Corollary
\ref{cor1.1}.

\begin{proof}
[Proof of Theorem \ref{thm1.1}]Let $h:K\rightarrow M$ be a Lipschitz
cubeulation. We may assume each cell in $K$ is a cube of unit size. Let
$\varepsilon_{M}>0$ be a small number such that%
\[
V_{2\varepsilon_{M}}\left(  M\right)  =\left\{  x\in\mathbb{R}^{l}:d\left(
x,N\right)  <2\varepsilon_{M}\right\}
\]
is a tubular neighborhood of $M$. Denote $\pi_{M}:V_{2\varepsilon_{M}}\left(
M\right)  \rightarrow M$ as the nearest point projection. For $\xi\in
B_{\varepsilon_{M}}^{l}$, we let $h_{\xi}\left(  x\right)  =\pi_{M}\left(
h\left(  x\right)  +\xi\right)  $ for $x\in\left\vert K\right\vert $, the
polytope of $K$. We may assume $\varepsilon_{M}$ is small enough such that all
$h_{\xi}$ are bi-Lipschitz maps. Replacing $h$ by $h_{\xi}$ when necessary, we
may assume $f=u\circ h\in\mathcal{W}^{1,2}\left(  K,N\right)  $. Then we may
find a $g\in C\left(  \left\vert K\right\vert ,N\right)  \cap\mathcal{W}%
^{1,2}\left(  K,N\right)  $ such that $\left[  g\circ h^{-1}\right]  =\alpha$
and $\left.  g\right\vert _{\left\vert K^{1}\right\vert }=\left.  f\right\vert
_{\left\vert K^{1}\right\vert }$ (see the proof of theorem 5.5 and theorem 6.1
in \cite{HL3}). For each cell $\Delta\in K$, let $y_{\Delta}$ be the center of
$\Delta$. For $x\in\Delta$, let $\left\vert x\right\vert _{\Delta}$ be the
Minkowski norm with respect to $y_{\Delta}$, that is%
\[
\left\vert x\right\vert _{\Delta}=\inf\left\{  t>0:y_{\Delta}+\frac
{x-y_{\Delta}}{t}\in\Delta\right\}  .
\]

\noindent\textbf{Step 1: }For every $\Delta\in K^{2}\backslash K^{1}$, we may
find a sequence $\phi_{i}\in C\left(  \Delta,N\right)  \cap W^{1,2}\left(
\Delta,N\right)  $ such that $\left.  \phi_{i}\right\vert _{\partial\Delta
}=\left.  g\right\vert _{\partial\Delta}$, $\phi_{i}\rightarrow\left.
f\right\vert _{\Delta}$ in $W^{1,2}\left(  \Delta,N\right)  $ and $d\phi
_{i}\rightarrow d\left(  \left.  f\right\vert _{\Delta}\right)  $ a.e. (see
lemma 4.4 in \cite{HL2}). For $x\in\Delta$, let%
\[
f_{i}\left(  x\right)  =\left\{
\begin{array}
[c]{cc}%
\phi_{i}\left(  x\right)  , & \left\vert x\right\vert _{\Delta}\geq\frac
{1}{2^{i}};\\
\phi_{i}\left(  y_{\Delta}+\frac{1}{2^{2i}\left\vert x\right\vert _{\Delta}%
}\frac{x-y_{\Delta}}{\left\vert x\right\vert _{\Delta}}\right)  , & \frac
{1}{2^{2i}}\leq\left\vert x\right\vert _{\Delta}\leq\frac{1}{2^{i}};\\
g\left(  y_{\Delta}+2^{2i}\left(  x-y_{\Delta}\right)  \right)  , & \left\vert
x\right\vert _{\Delta}\leq\frac{1}{2^{2i}}.
\end{array}
\right.
\]
It is clear that $f_{i}\rightharpoonup\left.  f\right\vert _{\Delta}$ in
$W^{1,2}\left(  \Delta,N\right)  $, $df_{i}\rightarrow d\left(  \left.
f\right\vert _{\Delta}\right)  $ a.e. on $\Delta$,%
\[
\left\vert df_{i}\right\vert _{L^{2}\left(  \Delta\right)  }\leq c\cdot\left(
\left\vert d\phi_{i}\right\vert _{L^{2}\left(  \Delta\right)  }+\left\vert
d\left(  \left.  g\right\vert _{\Delta}\right)  \right\vert _{L^{2}\left(
\Delta\right)  }\right)  \leq c\left(  f,g\right)
\]
and $f_{i}\in C\left(  \left\vert K^{2}\right\vert ,N\right)  $. In addition,
if we define $h_{2,i}:\Delta\times\left[  0,1\right]  \rightarrow N$ by%
\[
h_{2,i}\left(  x,t\right)  =\left\{
\begin{array}
[c]{cc}%
\phi_{i}\left(  x\right)  , & \left\vert x\right\vert _{\Delta}\geq\frac
{1}{2^{i}}+\frac{2^{i}-1}{2^{i}}t;\\
\phi_{i}\left(  y_{\Delta}+\frac{\left(  \frac{1}{2^{i}}+\frac{2^{i}-1}{2^{i}%
}t\right)  ^{2}}{\left\vert x\right\vert _{\Delta}}\frac{x-y_{\Delta}%
}{\left\vert x\right\vert _{\Delta}}\right)  , & \left(  \frac{1}{2^{i}}%
+\frac{2^{i}-1}{2^{i}}t\right)  ^{2}\leq\left\vert x\right\vert _{\Delta}%
\leq\frac{1}{2^{i}}+\frac{2^{i}-1}{2^{i}}t;\\
g\left(  y_{\Delta}+\frac{x-y_{\Delta}}{\left(  \frac{1}{2^{i}}+\frac{2^{i}%
-1}{2^{i}}t\right)  ^{2}}\right)  , & \left\vert x\right\vert _{\Delta}%
\leq\left(  \frac{1}{2^{i}}+\frac{2^{i}-1}{2^{i}}t\right)  ^{2}.
\end{array}
\right.
\]
Then by Lemma \ref{lem2.2}, we know $h_{2,i}\in W^{1,2}\left(  \Delta
\times\left[  0,1\right]  ,N\right)  $,
\[
\left\vert dh_{2,i}\right\vert _{L^{2}\left(  \Delta\times\left[  0,1\right]
\right)  }\leq c\cdot\left(  \left\vert d\phi_{i}\right\vert _{L^{2}\left(
\Delta\right)  }+\left\vert d\left(  \left.  g\right\vert _{\Delta}\right)
\right\vert _{L^{2}\left(  \Delta\right)  }\right)  \leq c\left(  f,g\right)
\]
and\ $h_{2,i}\in C\left(  \left\vert K^{2}\right\vert \times\left[
0,1\right]  ,N\right)  $.

\noindent\textbf{Step 2: }Assume for some $2\leq k\leq n-1$, we have a
sequence $f_{i}\in C\left(  \left\vert K^{k}\right\vert ,N\right)
\cap\mathcal{W}^{1,2}\left(  K^{k},N\right)  $ and $h_{k,i}\in C\left(
\left\vert K^{k}\right\vert \times\left[  0,1\right]  ,N\right)  $ such that
for each $\Delta\in K^{k}$, $f_{i}\rightharpoonup\left.  f\right\vert
_{\Delta}$ in $W^{1,2}\left(  \Delta,N\right)  $, $h_{k,i}\in W^{1,2}\left(
\Delta\times\left[  0,1\right]  ,N\right)  $,%
\begin{equation}
\left\vert d\left(  \left.  f_{i}\right\vert _{\Delta}\right)  \right\vert
_{L^{2}\left(  \Delta\right)  }\leq c\left(  f,g\right)  ,\quad\left\vert
dh_{k,i}\right\vert _{L^{2}\left(  \Delta\times\left[  0,1\right]  \right)
}\leq c\left(  f,g\right) \label{eq3.1}%
\end{equation}
and $h_{k,i}\left(  x,0\right)  =f_{i}\left(  x\right)  $, $h_{k,i}\left(
x,1\right)  =g\left(  x\right)  $ for $x\in\left\vert K^{k}\right\vert $.
Since for every $\Delta\in K^{k+1}\backslash K^{k}$, $f_{i}\rightharpoonup
\left.  f\right\vert _{\partial\Delta}$ in $W^{1,2}\left(  \partial
\Delta,N\right)  $, for fixed $j$ by Lemma \ref{lem2.1} we may find a
$n_{j}\geq j$ such that for each $\Delta\in K^{k+1}\backslash K^{k}$, there
exists a $w_{j}\in W^{1,2}\left(  \partial\Delta\times\left[  0,2^{-j}\right]
,N\right)  $ with $w_{j}\left(  x,0\right)  =f\left(  x\right)  $,
$w_{j}\left(  x,\frac{1}{2^{j}}\right)  =f_{n_{j}}\left(  x\right)  $ and%
\[
\left\vert dw_{j}\right\vert _{L^{2}\left(  \partial\Delta\times\left(
0,\frac{1}{2^{j}}\right)  \right)  }\leq\frac{c\left(  n\right)  }{2^{\frac
{j}{2}}}\left(  \left\vert d\left(  \left.  f\right\vert _{\partial\Delta
}\right)  \right\vert _{L^{2}\left(  \partial\Delta\right)  }+\left\vert
df_{n_{j}}\right\vert _{L^{2}\left(  \partial\Delta\right)  }+1\right)
\leq\frac{c\left(  f,g\right)  }{2^{\frac{j}{2}}}.
\]
Without loss of generality, we may replace $f_{i}$ by $f_{n_{i}}$ and
$h_{k,i}$ by $h_{k,n_{i}}$. Fix a $\Delta\in K^{k+1}\backslash K^{k}$. For
$x\in\Delta$, let%
\[
\psi_{i}\left(  x\right)  =\left\{
\begin{array}
[c]{cc}%
f\left(  y_{\Delta}+\frac{2^{i}\left(  x-y_{\Delta}\right)  }{2^{i}-1}\right)
, & \left\vert x\right\vert _{\Delta}\leq\frac{2^{i}-1}{2^{i}};\\
w_{i}\left(  y_{\Delta}+\frac{x-y_{\Delta}}{\left\vert x\right\vert _{\Delta}%
},\left\vert x\right\vert _{\Delta}-\frac{2^{i}-1}{2^{i}}\right)  , &
\frac{2^{i}-1}{2^{i}}\leq\left\vert x\right\vert _{\Delta}\leq1.
\end{array}
\right.
\]
Then $\left.  \psi_{i}\right\vert _{\left\vert K^{k}\right\vert }=f_{i}$ and
$\psi_{i}\rightarrow\left.  f\right\vert _{\Delta}$ in $W^{1,2}\left(
\Delta,N\right)  $ as $i\rightarrow\infty$ for each $\Delta\in K^{k+1}%
\backslash K^{k}$. By Theorem \ref{thmPR} and (\ref{eq3.1}) (use $h_{k,i}$ and
$g$ for the needed \textquotedblleft$v$\textquotedblright\ in Theorem
\ref{thmPR}, one may refer to lemma 9.8 of \cite{HL3}), for every $\Delta\in
K^{k+1}\backslash K^{k}$, we may find $\phi_{i}\in C\left(  \Delta,N\right)
\cap W^{1,2}\left(  \Delta,N\right)  $ such that $\left.  \phi_{i}\right\vert
_{\partial\Delta}=\left.  f_{i}\right\vert _{\partial\Delta}$, $\left\vert
\phi_{i}-\psi_{i}\right\vert _{L^{2}\left(  \Delta\right)  }<\frac{1}{2^{i}}$,
$\left\vert d\phi_{i}\right\vert _{L^{2}\left(  \Delta\right)  }\leq c\left(
f,g\right)  $ and
\[
\int_{M}\frac{\left\vert d\phi_{i}-d\psi_{i}\right\vert }{1+\left\vert
d\phi_{i}-d\psi_{i}\right\vert }d\mathcal{H}^{k+1}\leq\frac{1}{2^{i}}.
\]
After passing to subsequence, we may assume $d\phi_{i}\rightarrow d\left(
\left.  f\right\vert _{\Delta}\right)  $ a.e. on $\Delta$. Fix a $\Delta\in
K^{k+1}\backslash K^{k}$, for any $x\in\Delta$, define%
\begin{align*}
g_{k+1,i}\left(  x\right)   & =\left\{
\begin{array}
[c]{cc}%
h_{k,i}\left(  y_{\Delta}+\frac{x-y_{\Delta}}{\left\vert x\right\vert
_{\Delta}},1+2\left(  \frac{1}{2}-\left\vert x\right\vert _{\Delta}\right)
\right)  , & \frac{1}{2}\leq\left\vert x\right\vert _{\Delta}\leq1;\\
g\left(  y_{\Delta}+2\left(  x-y_{\Delta}\right)  \right)  , & \left\vert
x\right\vert _{\Delta}\leq\frac{1}{2},
\end{array}
\right. \\
f_{i}\left(  x\right)   & =\left\{
\begin{array}
[c]{cc}%
\phi_{i}\left(  x\right)  , & \left\vert x\right\vert _{\Delta}\geq\frac
{1}{2^{i}};\\
\phi_{i}\left(  y_{\Delta}+\frac{1}{2^{2i}\left\vert x\right\vert _{\Delta}%
}\frac{x-y_{\Delta}}{\left\vert x\right\vert _{\Delta}}\right)  , & \frac
{1}{2^{2i}}\leq\left\vert x\right\vert _{\Delta}\leq\frac{1}{2^{i}};\\
g_{k+1,i}\left(  y_{\Delta}+2^{2i}\left(  x-y_{\Delta}\right)  \right)  , &
\left\vert x\right\vert _{\Delta}\leq\frac{1}{2^{2i}},
\end{array}
\right. \\
\widetilde{h}_{k+1,i}\left(  x,t\right)   & =\left\{
\begin{array}
[c]{cc}%
\phi_{i}\left(  x\right)  , & \left\vert x\right\vert _{\Delta}\geq\frac
{1}{2^{i}}+\frac{2^{i}-1}{2^{i}}t;\\
\phi_{i}\left(  y_{\Delta}+\frac{\left(  \frac{1}{2^{i}}+\frac{2^{i}-1}{2^{i}%
}t\right)  ^{2}}{\left\vert x\right\vert _{\Delta}}\frac{x-y_{\Delta}%
}{\left\vert x\right\vert _{\Delta}}\right)  , & \left(  \frac{1}{2^{i}}%
+\frac{2^{i}-1}{2^{i}}t\right)  ^{2}\leq\left\vert x\right\vert _{\Delta}%
\leq\frac{1}{2^{i}}+\frac{2^{i}-1}{2^{i}}t;\\
g_{k+1,i}\left(  y_{\Delta}+\frac{x-y_{\Delta}}{\left(  \frac{1}{2^{i}}%
+\frac{2^{i}-1}{2^{i}}t\right)  ^{2}}\right)  , & \left\vert x\right\vert
_{\Delta}\leq\left(  \frac{1}{2^{i}}+\frac{2^{i}-1}{2^{i}}t\right)  ^{2},
\end{array}
\right. \\
\widetilde{\widetilde{h}}_{k+1,i}\left(  x,t\right)   & =\left\{
\begin{array}
[c]{cc}%
h_{k,i}\left(  y_{\Delta}+\frac{x-y_{\Delta}}{\left\vert x\right\vert
_{\Delta}},1+2\left(  \frac{1+t}{2}-\left\vert x\right\vert _{\Delta}\right)
\right)  , & \frac{1+t}{2}\leq\left\vert x\right\vert _{\Delta}\leq1;\\
g\left(  y_{\Delta}+\frac{2}{1+t}\left(  x-y_{\Delta}\right)  \right)  , &
\left\vert x\right\vert _{\Delta}\leq\frac{1+t}{2},
\end{array}
\right.
\end{align*}
and%
\[
h_{k+1,i}\left(  x,t\right)  =\left\{
\begin{array}
[c]{cc}%
\widetilde{h}_{k+1,i}\left(  x,2t\right)  , & 0\leq t\leq\frac{1}{2};\\
\widetilde{\widetilde{h}}_{k+1,i}\left(  x,2t-1\right)  , & \frac{1}{2}\leq
t\leq1.
\end{array}
\right.
\]
Simple calculations show that for any $\Delta\in K^{k+1}\backslash K^{k}$,
$f_{i}\rightharpoonup\left.  f\right\vert _{\Delta}$ in $W^{1,2}\left(
\Delta,N\right)  $, $df_{i}\rightarrow d\left(  \left.  f\right\vert _{\Delta
}\right)  $ a.e. on $\Delta$, $h_{k+1,i}\in W^{1,2}\left(  \Delta\times\left[
0,1\right]  ,N\right)  $,%
\[
\left\vert df_{i}\right\vert _{L^{2}\left(  \Delta\right)  }\leq c\left(
f,g\right)  ,\quad\left\vert dh_{k+1,i}\right\vert _{L^{2}\left(  \Delta
\times\left[  0,1\right]  \right)  }\leq c\left(  f,g\right)
\]
and $h_{k+1,i}\left(  x,0\right)  =f_{i}\left(  x\right)  $, $h_{k+1,i}\left(
x,1\right)  =g\left(  x\right)  $ for $x\in\left\vert K^{k+1}\right\vert $.
Hence we finish when we reach $f_{i}\in C\left(  \left\vert K\right\vert
,N\right)  \cap\mathcal{W}^{1,2}\left(  K,N\right)  $ and $h_{n,i}\in C\left(
\left\vert K\right\vert \times\left[  0,1\right]  ,N\right)  $. Let
$v_{i}=f_{i}\circ h^{-1}$. Then it is clear that $v_{i}\in C\left(
M,N\right)  \cap W^{1,2}\left(  M,N\right)  $, $\left[  v_{i}\right]  =\alpha
$, $\left\vert v_{i}-u\right\vert _{L^{2}\left(  M\right)  }\rightarrow0$,
$\left\vert dv_{i}\right\vert _{L^{2}\left(  M\right)  }\leq c\left(
u,g\right)  $ and $dv_{i}\rightarrow du$ a.e. on $M$. Hence, we may find
$u_{i}\in C^{\infty}\left(  M,N\right)  $ such that $\left\vert u_{i}%
-u\right\vert _{L^{2}\left(  M\right)  }\rightarrow0$, $\left\vert
du_{i}\right\vert _{L^{2}\left(  M\right)  }\leq c\left(  u,g\right)  $,
$\left[  u_{i}\right]  =\alpha$ and $du_{i}\rightarrow du$ a.e. on $M$. In
particular, this shows%
\[
H_{W}^{1,2}\left(  M,N\right)  \supset\left\{  u\in W^{1,2}\left(  M,N\right)
:u_{\#,2}\left(  h\right)  \text{ has a continuous extension to }M\text{
w.r.t. }N\right\}  .
\]
The other direction of inclusion was proved in section 7 of \cite{HL2}. To see%
\[
H_{W}^{1,2}\left(  M,N\right)  =\left\{  u\in W^{1,2}\left(  M,N\right)
:u\text{ may be connected to some smooth maps}\right\}  ,
\]
we only need to use the above proved equality and proposition 5.2 of
\cite{HL2}, which shows%
\begin{align*}
& \left\{  u\in W^{1,2}\left(  M,N\right)  :u_{\#,2}\left(  h\right)  \text{
has a continuous extension to }M\text{ w.r.t. }N\right\} \\
& =\left\{  u\in W^{1,2}\left(  M,N\right)  :u\text{ may be connected to some
smooth maps}\right\}  .
\end{align*}

\end{proof}

We remark that many constructions above are motivated from section 5 and
section 6 of \cite{HL3}.

\begin{proof}
[Proof of Corollary \ref{cor1.1}]This follows from Theorem \ref{thm1.1} and
corollary 5.4 of \cite{HL2}.
\end{proof}


\begin{thebibliography}{9}                                                                                                %
\bibitem {H}P. Hajlasz. Approximation of Sobolev mappings. \textit{Nonlinear
Anal} \textbf{22} (1994), no. 12, 1579--1591.

\bibitem {HL1}F. B. Hang and F. H. Lin. Topology of Sobolev mappings.
\textit{Math Res Lett} \textbf{8} (2001), no. 3, 321--330.

\bibitem {HL2}F. B. Hang and F. H. Lin. Topology of Sobolev mappings II.
\textit{Acta Math} \textbf{191} (2003), no. 1, 55--107.

\bibitem {HL3}F. B. Hang and F. H. Lin. Topology of Sobolev mappings III.
\textit{Comm Pure Appl Math} \textbf{56} (2003), no. 10, 1383--1415.

\bibitem {HR}R. Hardt and T. Riviere. Connecting topological Hopf
singularities. \textit{Annali Sc Norm Sup Pisa}, \textbf{2} (2003), no. 2, 287--344.

\bibitem {L}S. Luckhaus. Partial Holder continuity for minima of certain
energies among maps into a Riemannian manifold. \textit{Indiana Univ Math J}
\textbf{37} (1988), 349--367.

\bibitem {PR}M. R. Pakzad and T. Riviere. Weak density of smooth maps for the
Dirichlet energy between manifolds. \textit{Geom Func Anal} \textbf{13}
(2003), no. 1, 223--257.

\bibitem {W}B. White. Homotopy classes in Sobolev spaces and the existence of
energy minimizing maps. \textit{Acta Math} \textbf{160} (1988), no. 1--2,
1--17.\textit{\ }
\end{thebibliography}
\end{document}